\documentclass[a4paper,11pt]{article}
\usepackage{amssymb}
\usepackage{amsmath}
\usepackage{latexsym,a4wide}
\usepackage{xy}

\input xypic
\input xy
\xyoption{all}

\newtheorem{example}{Example}[section]
\newtheorem{defn}[example]{Definition}
\newtheorem{prop}[example]{Proposition}
\newtheorem{thm}[example]{Theorem}

\newenvironment{pf}{\noindent \textbf{Proof:}}{\rule{0em}{1ex}\hfill$\Box$\mbox{}}

\def\leq{\leqslant}

\input{tcilatex}

\date{}
\begin{document}

\author{Z. Arvasi and E. Ulualan}
\title{Homotopical Aspects of Commutative
Algebras}
\maketitle

\begin{abstract}
This article investigates the homotopy theory of simplicial
commutative algebras with a  view to homological applications.
\footnote{Quadratic Modules, Simplicial Commutative Algebras,
2-Crossed Modules, Crossed Squares.}

\end{abstract}

\section*{Introduction}

The original motivation for this article was to see what parts of
the group theoretic case of crossed homotopical algebra
generalised to the context of commutative algebras and to see how
existing parts of commutative algebra might interact with the
analogue. The hope was for a clarification of the group theoretic
situation as well as perhaps introducing `new' tools into
commutative algebra.

Simplicial commutative algebras occupy a place somewhere between
homological algebra, homotopy theory, algebraic K-theory and
algebraic geometry. In each sector they have played a significant
part in developments over quite a lengthy period of time. Their
own internal structure has however been studied in \cite{Arvasi,
Arvasi1, Arvasi2}. The present work gives some lights on the
3-types of simplicial commutative algebras and will apply the
results in various mainly homological settings.

Crossed modules of groups were introduced by Whitehead \cite{wayted}. They
model homotopy types with trivial homotopy groups in dimensions bigger than
2. Algebraic models for connected 3-types in terms of crossed module of
groups of length 2 are 2-crossed module defined by Conduch{\'{e}} \cite{con}
and crossed square defined by Guin-Walery, Loday \cite{walery}. The
commutative algebra version of these structures are respectively defined by
Grandje\'{a}n-Vale \cite{vale} and Ellis \cite{elis}. The first author and
Porter (cf. \cite{Arvasi, Arvasi1, Arvasi2}) also gave the relations between
these constructions and simplicial algebras. For an alternative model, the
notion of quadratic module of groups is defined by Baues (cf. \cite{baus1}).
This is a 2-crossed module with additional nilpotent conditions.

In this paper, we define the commutative algebra version of quadratic
modules and construct a functor from the category of simplicial commutative
algebras to that of quadratic modules by using higher dimensional Peiffer
elements. In this construction we see the role of the hypercrossed complex
pairings. Furthermore, we give a link between quadratic modules, 2-crossed
modules, crossed squares and simplicial commutative algebras with Moore
complex of length 2. Another algebraic model of 3-types is `braided regular
crossed module' introduced by Brown and Gilbert in \cite{rb-ng}. This model
will be analyzed in a separate paper.

Quillen \cite{q} and Illusie \cite{ill} both discuss the basic homotopical
algebra of simplicial algebras and their application in deformation theory.
Andr\'{e} \cite{andre} gives a detailed examination of their construction
and applies them to cohomology via the cotangent complex construction.

Thus the situations which are examined in this paper and related models such
as \cite{Arvasi}, \cite{Arvasi1} can be summarised in the following diagram
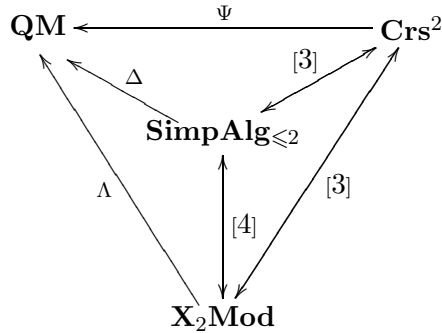
\begin{equation*}
\xymatrix{\textbf{QM}&&\textbf{Crs$^2$}\ar[dddl]^{\cite{Arvasi}}\ar@<0.01ex>[ll]_{\Psi}\ar[dl]_{%
\cite{Arvasi}}\\
&\textbf{SimpAlg$_{\leq
2}$}\ar[ur]\ar[dd]^{\cite{Arvasi1}}\ar@<0.99ex>[ul]_{\Delta}&\\ \\
&\textbf{X$_2$Mod}\ar[uu]\ar@<0.99ex>[uuul]^{\Lambda}\ar[uuur]&&}
\end{equation*}%
where the functors $\Lambda $, $\Delta ,$ $\Psi $ are given by Propositions %
\ref{quad}, \ref{simp} and \ref{crs} respectively and the numbers in the
diagram correspond to the references.

\section{\label{pre} Preliminaries}

In what follows `algebras' will be commutative algebras over an unspecified
commutative ring, \textbf{k}, but for convenience are not required to have a
multiplicative identity. The category of commutative algebras will be
denoted by $\mathbf{Alg}$.

\textbf{Simplicial Commutative Algebras}

A simplicial (commutative) algebra $\mathbf{E}$ consists of a family of
algebras $\{E_{n}\}$ together with face and degeneracy maps $%
d_{i}=d_{i}^{n}:E_{n}\rightarrow E_{n-1},$ \ \ $0\leqslant i\leqslant n$, \ $%
(n\neq 0)$ and $s_{i}=s_{i}^{n}:E_{n}\rightarrow E_{n+1}$, \ \ $0\leqslant
i\leqslant n$, satisfying the usual simplicial identities given in Andr\'{e}
\cite{andre} or Illusie \cite{ill} for example. It can be completely
described as a functor $\mathbf{E}:\Delta ^{op}\rightarrow \mathbf{Alg}$
where $\Delta $ is the category of finite ordinals $[n]=\{0<1<\cdots <n\}$
and increasing maps. We denote the category of simplicial algebras by $%
\mathbf{SimpAlg.}$

Given a simplicial algebra $\mathbf{E}$, the Moore complex $(\mathbf{NE}%
,\partial )$ of $\mathbf{E}$ is the chain complex defined by
\begin{equation*}
{NE}_{n}=\ker d_{0}^{n}\cap \ker d_{1}^{n}\cap \cdots \cap \ker d_{n-1}^{n}
\end{equation*}%
with $\partial _{n}:NE_{n}\rightarrow NE_{n-1}$ induced from $d_{n}^{n}$ by
restriction.

The $n$th \emph{homotopy module} $\pi _{n}(\mathbf{E})$ of \textbf{E} is the
$n$th homology of the Moore complex of $\mathbf{E}$, i.e.
\begin{equation*}
\pi _{n}\mathbf{(E)}\cong H_{n}(\mathbf{NE},\partial
)=\bigcap\limits_{i=0}^{n}\ker d_{i}^{n}/d_{n+1}^{n+1}\left(
\bigcap\limits_{i=0}^{n}\ker d_{i}^{n+1}\right)
\end{equation*}%
\newline
We say that the Moore complex $\mathbf{NE}$ of a simplicial algebra is of
length $k$ if $NE_{n}=0$ for all $n\geqslant k+1,$ so that a Moore complex
of length $k$ is also of length $l$ for $l\geqslant k.$ We denote thus the
category of simplicial algebras with Moore complex of length $k$ by $\mathbf{%
SimpAlg}_{\leqslant k}$.


The following terminologies and notations are derived from \cite{c} and the
published version, \cite{cc}, of the analogous group theoretic case. For
detailed investigation see \cite{Arvasi1}.

For the ordered set $[n]=\{0<1<\cdots <n\}$, let $\alpha
_{i}^{n}:[n+1]\rightarrow \lbrack n]$ be the increasing surjective map given
by
\begin{equation*}
\alpha _{i}^{n}(j)=\left\{
\begin{array}{ll}
j & \text{if }j\leqslant i, \\
j-1 & \text{if }j>i.%
\end{array}%
\right.
\end{equation*}%
Let $S(n,n-r)$ be the set of all monotone increasing surjective maps from $%
[n]$ to $[n-r]$. This can be generated from the various $\alpha _{i}^{n}$ by
composition. The composition of these generating maps is subject to the
following rule: $\alpha _{j}\alpha _{i}=\alpha _{i-1}\alpha _{j},$ $j<i$.
This implies that every element $\alpha \in S(n,n-r)$ has a unique
expression as $\alpha =\alpha _{i_{1}}\circ \alpha _{i_{2}}\circ \cdots
\circ \alpha _{i_{r}}$ with $0\leqslant i_{1}<i_{2}<\cdots <i_{r}\leqslant
n-1$, where the indices $i_{k}$ are the elements of $[n]$ such that $%
\{i_{1},\ldots ,i_{r}\}=\{i:\alpha (i)=\alpha (i+1)\}$. We thus can identify
$S(n,n-r)$ with the set $\{(i_{r},\ldots ,i_{1}):0\leqslant
i_{1}<i_{2}<\cdots <i_{r}\leqslant n-1\}$. In particular, the single element
of $S(n,n)$, defined by the identity map on $[n]$, corresponds to the empty
0-tuple ( ) denoted by $\emptyset _{n}$. Similarly the only element of $%
S(n,0)$ is $(n-1,n-2,\ldots ,0)$. For all $n\geqslant 0$, let
\begin{equation*}
S(n)=\bigcup_{0\leqslant r\leqslant n}S(n,n-r).
\end{equation*}%
We say that $\alpha =(i_{r},\ldots ,i_{1})<\beta =(j_{s},\ldots ,j_{1})$ in $%
S(n),$ if $i_{1}=j_{1},\ldots ,i_{k}=j_{k}$ but $i_{k+1}>j_{k+1},(k\geqslant
0)$ or if $i_{1}=j_{1},\ldots ,i_{r}=j_{r}$ and $r<s$. \newline
This makes $S(n)$ an ordered set.

\textbf{Hypercrossed Complex Pairings}

We give the following statements from \cite{Arvasi1}. For details see \cite%
{cc} (group case) and \cite{Arvasi1}. We define a set $P(n)$ consisting of
pairs of elements $(\alpha ,\beta )$ from $S(n)$ with $\alpha \cap \beta
=\emptyset $ and $\beta <\alpha $ where $\alpha =(i_{r},\ldots ,i_{1})$, $%
\beta =(j_{s},\ldots ,j_{1})\in S(n)$. The \textbf{k}-linear morphisms that
we will need,
\begin{equation*}
\{C_{\alpha ,\beta }:NE_{n-\#\alpha }\otimes NE_{n-\#\beta }\rightarrow
NE_{n}:(\alpha ,\beta )\in P(n),n\geqslant 0\}
\end{equation*}%
are given as composites:
\begin{eqnarray*}
C_{\alpha ,\beta }(x_{\alpha }\otimes y_{\beta }) &=&p\mu (s_{\alpha
}\otimes s_{\beta })(x_{\alpha }\otimes y_{\beta }) \\
&=&p(s_{\alpha }(x_{\alpha })s_{\beta }(x_{\beta })) \\
&=&(1-s_{n-1}d_{n-1})\cdots (1-s_{0}d_{0})(s_{\alpha }(x_{\alpha })s_{\beta
}(x_{\beta })),
\end{eqnarray*}%
where $s_{\alpha }=s_{i_{r}}\ldots s_{i_{1}}:NE_{n-\#\alpha }\rightarrow
E_{n}$, $s_{\beta }=s_{j_{s}}\ldots s_{j_{1}}:NE_{n-\#\beta }\rightarrow
E_{n},$ $p:E_{n}\rightarrow NE_{n}$ is defined by composite projections $%
p=p_{n-1}\ldots p_{0}$ with
\begin{equation*}
p_{j}=1-s_{j}d_{j}\text{ for }j=0,1,\ldots ,n-1
\end{equation*}%
and $\mu :E_{n}\otimes E_{n}\rightarrow E_{n}$ denotes multiplication.

We will now consider that the ideal $I_{n}$ in $E_{n}$ such that generated
by all elements of the form;
\begin{equation*}
C_{\alpha ,\beta }(x_{\alpha }\otimes y_{\beta })
\end{equation*}
where $x_{\alpha }\in NE_{n-\#\alpha }$ and $y_{\beta }\in NE_{n-\#\beta }$
and for all $(\alpha ,\beta )\in P(n)$.

\begin{prop}
\emph{(\cite{Arvasi1})} Let \textbf{E} be simplicial algebra and $n>0$, and $%
D_{n}$ the ideal in $E_{n}$ generated by degenerate elements. We suppose $%
E_{n}=D_{n}$, and let $I_{n}$ be the ideal generated by elements of the form
\begin{equation*}
C_{\alpha ,\beta }(x_{\alpha }\otimes y_{\beta })\qquad \text{with }(\alpha
,\beta )\in P(n)
\end{equation*}
where $x_{\alpha }\in NE_{n-\#\alpha },y_{\beta }\in NE_{n-\#\beta }$ with $%
1\leqslant r,s\leqslant n.$ Then,
\begin{equation*}
\partial _{n}(NE_{n})=\partial _{n}(I_{n}).
\end{equation*}
\end{prop}

Now according to above proposition for $n=2$, $3,$ we show the image of $%
I_{n}$ by $\partial _{n}$ what it looks like.

We suppose $D_{2}=E_{2}$. We take $\beta =(1),\alpha =(0)$ and $x,y\in
NE_{1} $ = $\ker d_{0}$. We know that the ideal $I_{2}$ is generated by
elements of the form
\begin{equation*}
C_{(1)(0)}(x\otimes y)=p_{1}p_{0}(s_{1}xs_{0}y)=s_{1}x(s_{1}y-s_{0}y).
\end{equation*}%
Then the image of $I_{2}$ by $\partial _{2}$ is $d_{2}[C_{(1)(0)}(x\otimes
y)]=x(y-s_{0}d_{1}y)$ where $x\in $ $\ker d_{0}$ and $y-s_{0}d_{1}y\in \ker
d_{1}$. Therefore the image of $I_{2}$ by $\partial _{2}$ is $\ker d_{0}\ker
d_{1}.$\newline
For $n=3,$ the linear morphisms are
\begin{equation*}
\begin{array}{lllll}
C_{(1,0),(2)}, &  & C_{(2,0),(1)}, &  & C_{(2,1),(0)} \\
C_{(2),(0)}, &  & C_{(2),(1)}, &  & C_{(1),(0)}.%
\end{array}%
\end{equation*}%
Then the ideal $I_{3}$ is generated by the elements; for $x\in NE_{1}$ and $%
y\in NE_{2}$
\begin{align*}
C_{(1,0),(2)}(x\otimes y) & = (s_{1}s_{0}x-s_{2}s_{0}x)s_{2}y \\
C_{(2,0),(1)}(x\otimes y) & = (s_{2}s_{0}x-s_{2}s_{1}x)(s_{1}y-s_{2}y) \\
C_{(2,0),(1)}(x\otimes y) & = s_{2}s_{1}x(s_{0}y-s_{1}y+s_{2}y)
\end{align*}
and for $x,y\in NE_{2}$
\begin{align*}
C_{(1),(0)}(x\otimes y) & = s_{1}x(s_{0}y-s_{1}y)+s_{2}xy \\
C_{(2),(0)}(x\otimes y) & = s_{2}xs_{0}y \\
C_{(1),(0)}(x\otimes y) & = s_{2}x(s_{1}y-s_{2}y).
\end{align*}%
Thus the image of $I_{3}$ by $\partial _{3}$ is
\begin{align*}
\partial _{3}(C_{(1,0)(2)}(x\otimes y)) & = (s_{1}s_{0}d_{1}x-s_{0}x)y\in
\ker d_{2}(\ker d_{0}\cap \ker d_{1}) \\
\partial _{3}(C_{(2,0)(1)}(x\otimes y)) & =
(s_{0}x-s_{1}x)(s_{1}d_{2}y-y)\in \ker d_{1}(\ker d_{0}\cap \ker d_{2}) \\
\partial _{3}(C_{(2,1)(0)}(x\otimes y)) & =
s_{1}x(s_{0}d_{2}y-s_{1}d_{2}y+y)\in \ker d_{0}(\ker d_{1}\cap \ker d_{2}) \\
\partial _{3}(C_{(2)(1)}(x\otimes y)) & = x(s_{1}d_{2}y-y)\in (\ker
d_{0}\cap \ker d_{1})(\ker d_{0}\cap \ker d_{2}) \\
\partial _{3}(C_{(2)(0)}(x\otimes y)) & = xs_{0}d_{2}y\in
K_{\{0,1\}}K_{\{1,2\}}+K_{\{0,1\}}K_{\{0,2\}} \\
\partial _{3}(C_{(1)(0)}(x\otimes y)) & =
s_{1}d_{2}xs_{0}d_{2}y-s_{1}d_{2}xs_{1}d_{2}y+xy \\
& \in K_{\{0,2\}}K_{\{1,2\}}+K_{\{0,1\}}K_{\{1,2\}}+K_{\{0,1\}}K_{\{0,2\}}.
\end{align*}

If $n=4$ then the image of the Moore complex of the simplicial algebra $%
\mathbf{E}$\textbf{\ }can be given in the form
\begin{equation*}
\partial _{n}(NE_{n})=\underset{I,J}{\sum }K_{I}K_{J},
\end{equation*}%
where $\emptyset \neq I,J\subset \lbrack n-1]=\{0,1,\dots ,n-2\}$ with $%
I\cup J=[n-1],$ and where
\begin{equation*}
K_{I}=\underset{i\in I}{\bigcap }\ker d_{i}\text{ and }K_{J}=\underset{J\in J%
}{\bigcap }\ker d_{j}.
\end{equation*}%
In general for $n>4$, there is an inclusion
\begin{equation*}
\sum_{I,J}K_{I}K_{J}\subset \partial _{n}(NE_{n}).
\end{equation*}

\section{\label{simp-2crs}2-Crossed Modules from Simplicial Algebras}

Crossed modules techniques give a very efficient way of handling information
about on a homotopy type. They correspond to a 2-type. The commutative
algebra analogue of crossed modules has been given by Porter in \cite%
{porter4}. Throughout this paper we denote an action of $r\in R$ on $c\in C$
by $c\cdot r$.

A \textit{crossed module } is an $R$-algebra homomorphism $\partial
:C\rightarrow R$ with the action of $R$ on $C$ such that $\partial (c\cdot
r)=\partial (c)r$ and $c^{\prime }\cdot \partial (c)=c^{\prime }c$ for all $%
r\in R$ and $c,c^{\prime }\in C.$ The second condition is called the \textit{%
Peiffer identity}. We will denote such a crossed module by $(C,R,\partial ).$
A morphism of crossed modules from $(C,R,\partial )$ to $(C^{\prime
},R^{\prime },\partial ^{\prime })$ is a pair of $\mathbf{k}$-algebra
homomorphisms, $\phi :C\rightarrow C^{\prime }$ and $\varphi :R\rightarrow
R^{\prime }$ such that $\phi (c\cdot r)=\phi (c)\varphi (r).$ We thus define
the category of crossed modules of commutative algebras denoting it by $%
\mathbf{XMod.}$

As we mentioned in introduction, the notion 2-crossed modules of groups was
introduced by Conduch\'{e} in \cite{con} as a model for connected 3-types.
He showed that the category of 2-crossed modules is equivalent to the
category of simplicial groups with Moore complex of length 2. Grandje\'{a}n
and Vale, \cite{vale}, gave the notion of 2-crossed modules on commutative
algebras. Now, we recall from \cite{vale} the definition of a 2-crossed
module:

A \textit{2-crossed module} of $\mathbf{k}$-algebras consists of a complex
of $C_{0}$-algebras
\begin{equation*}
\xymatrix{C_2 \ar[r]^-{\partial _{2}}&C_1 \ar[r]^-{\partial _{1}}&C_0}
\end{equation*}%
and $\partial _{2},$ $\partial _{1}$ morphisms of $C_{0}$-algebras, where
the algebra $C_{0}$ acts on itself by multiplication such that
\begin{equation*}
\partial _{2}:C_{2}\longrightarrow C_{1}
\end{equation*}%
is a crossed module. Thus $C_{1}$ acts on $C_{2}$ via $C_{0}$ and we require
that for all $x\in C_{2},$ $y\in C_{1},$ $z\in C_{0}$ that $(xy)z=x(yz).$
Further there is a $C_{0}$-bilinear function giving
\begin{equation*}
\{-\otimes -\}:C_{1}\otimes _{C_{0}}C_{1}\longrightarrow C_{2}
\end{equation*}%
called a Peiffer lifting, which satisfies the following axioms
\begin{equation*}
\begin{array}{crcl}
\mathbf{2CM1)} & \partial _{2}\{y_{0}\otimes y_{1}\} & = &
y_{0}y_{1}-y_{0}\cdot \partial _{1}y_{1} \\
\mathbf{2CM2)} & \{\partial _{2}(x_{1})\otimes \partial _{2}(x_{2})\} & = &
x_{1}x_{2} \\
\mathbf{2CM3)} & \{y_{0}\otimes y_{1}y_{2}\} & = & \{y_{0}y_{1}\otimes
y_{2}\}+\{y_{0}\otimes y_{1}\}\cdot \partial _{1}y_{2} \\
\mathbf{2CM4)} & a)\{\partial _{2}x\otimes y\} & = & x\cdot y-x\cdot
\partial _{1}y \\
& b)\{y\otimes \partial _{2}x\} & = & x\cdot y \\
\mathbf{2CM5)} & \{y_{0}\otimes y_{1}\}\cdot z & = & \{y_{0}\cdot z\otimes
y_{1}\}=\{y_{0}\otimes y_{1}\cdot z\}%
\end{array}%
\end{equation*}%
\newline
for all $x,x_{1},x_{2}\in C_{2}$, $y,y_{0},y_{1},y_{2}\in C_{1}$ and $z\in
C_{0}.$

A morphism of 2-crossed modules of algebras may be pictured by the diagram
\begin{equation*}
\xymatrix{ C_{2} \ar[d]_{f_{2}} \ar[r]^{\partial_{2}} & C_{1} \ar[d]_{f_{1}}
\ar[r]^{\partial_{1}} & C_{0}\ar[d]_{f_{0}} \\ C_{2}^{\prime}
\ar[r]_{\partial^{\prime}_{2}} & C_{1}^{\prime}
\ar[r]_{\partial^{\prime}_{1}} & C_{0}^{\prime} }
\end{equation*}
such that
\begin{equation*}
f_{0}\partial _{1}=\partial _{1}^{\prime }f_{1},\quad f_{1}\partial
_{2}=\partial _{2}^{\prime }f_{2}
\end{equation*}
and such that
\begin{equation*}
f_{1}(c_{0}\cdot c_{1})=f_{0}(c_{0})\cdot f_{1}(c_{1}),\quad
f_{2}(c_{0}\cdot c_{2})=f_{0}(c_{0})\cdot f_{2}(c_{2})
\end{equation*}
and
\begin{equation*}
\{\otimes \}f_{1}\otimes f_{1}=f_{2}\{\otimes \}
\end{equation*}
for all $c_{2}\in C_{2},$ $c_{1}\in C_{1},$ $c_{0}\in C_{0}$ .

We denote the category of 2-crossed module by $\mathbf{X}_{2}\mathbf{Mod}$%
\textbf{.}

In \cite{Arvasi1}, the first author with Porter studied the truncated
simplicial algebras and saw what properties that has. Later, they turned to
a simplicial algebra $\mathbf{E}$ which is 2-truncated, i.e., its Moore
complex look like;%
\begin{equation*}
\xymatrix{ \cdots \ar[r]&0\ar[r]&NE_2 \ar[r]^{\partial_2} & NE_{1}
\ar[r]^{\partial_{1}} &NE_{0}}
\end{equation*}%
and they then showed the following result:

\begin{thm}
The category $\mathbf{X}_{2}\mathbf{Mod}$ of 2-crossed modules is equivalent
to the category $\mathbf{SimpAlg}_{\leqslant 2}$ of simplicial algebras with
Moore complex of length 2.
\end{thm}

\section{\label{square-2crs}2-Crossed Modules from Crossed Squares}

Crossed squares were initially defined by Guin-Wal\'{e}ry and Loday in \cite%
{walery}. The commutative algebra analogue of crossed squares has been
studied by Ellis (cf. \cite{elis}).

Now, we will recall the definition (due to Ellis \cite{elis}) of a crossed
square of algebras.\newline
A\textit{\ crossed square} is a commutative diagram
\begin{equation*}
\begin{array}{c}
\xymatrix{ L\ar[r]^{\lambda}\ar[d]_{\lambda'}& M\ar[d]^{\mu}\\ N
\ar[r]_{\nu}&R}%
\end{array}%
\end{equation*}%
together with the actions of $R$ on $L,M$ and $N.$ There are thus
commutative actions of $M$ on $L$ and $N$ via $\mu $ and a function
\begin{equation*}
h:M\times N\rightarrow L
\end{equation*}%
the $h$-map. This data must satisfy the following axioms: \newline
$1.\qquad $The $\lambda ,\lambda ^{\prime },\mu ,\nu $ and $\mu \lambda =\nu
\lambda ^{\prime }$ are crossed modules;\newline
$2.\qquad $The maps $\lambda ,\lambda ^{\prime }$ preserve the action of $R$;%
\newline
$3.\qquad kh(m,n)=h(km,n)=h(m,kn);$\newline
$4.\qquad h(m+m^{\prime },n)=h(m,n)+h(m^{\prime },n);$\newline
$5.\qquad h(m,n+n^{\prime })=h(m,n)+h(m,n^{\prime });$\newline
$6.\qquad r\cdot h(m,n)=h(r\cdot m,n)=h(m,r\cdot n);$\newline
$7.\qquad \lambda h(m,n)=m\cdot n;$\newline
$8.\qquad \lambda ^{\prime }h(m,n)=n\cdot m;$\newline
$9.\qquad h(m,\lambda ^{\prime }l)=m\cdot l;$\newline
$10.\qquad h(\lambda l,n)=n\cdot l$ \newline
for all $l\in L,m,m^{\prime }\in M$ and $n,n^{\prime }\in N,$ $r\in R$, $%
k\in \mathbf{k}$

We denote such a crossed square by $(L,M,N,R)$. A morphism of crossed
squares $\Phi :(L,M,N,R)\rightarrow (L^{\prime },M^{\prime },N^{\prime
},R^{\prime })$ consists of homomorphisms
\begin{equation*}
\begin{array}{cccc}
\Phi _{L}: & L & \longrightarrow & L^{\prime } \\
\Phi _{N}: & N & \longrightarrow & N^{\prime }%
\end{array}
\quad
\begin{array}{cccc}
\Phi _{M}: & M & \longrightarrow & M^{\prime } \\
\Phi _{R}: & R & \longrightarrow & R^{\prime }%
\end{array}%
\end{equation*}
such that the cube of homomorphisms is commutative
\begin{equation*}
\Phi _{L}h(m,n)=h(\Phi _{M}m,\Phi _{N}n)
\end{equation*}
with $m\in M$ and $n\in N,$ and the homomorphisms $\Phi _{L},\Phi _{M},\Phi
_{N}$ are $\Phi _{R}$-equivariant. The category of crossed squares will be
denoted by $\mathbf{Crs}^{2}.$

Conduch\'{e} (in a private communication with Brown and (see also published
version \cite{Con1})) gives a construction of a 2-crossed module from a
crossed square of groups. On the other hand, the first author gave a neat
description of the passage from a crossed square to a 2-crossed module of
algebras in \cite{Arvasi}. He constructed a 2-crossed module from a crossed
square of commutative algebras
\begin{equation*}
\xymatrix{ L \ar[d]_{\lambda'} \ar[r]^{\lambda} & M \ar[d]^{\mu} \\ N
\ar[r]_{\nu} & R }
\end{equation*}%
as
\begin{equation*}
\xymatrix{L\ar[rr]^-{(-\lambda,\lambda ^{\prime })}&&M\rtimes
N\ar[rr]^-{\mu+\nu}&&R }
\end{equation*}%
analogue to that given by Conduch\'{e} in the group case (cf. \cite{Con1}).
This construction can be briefly summarised as follows:

Apply the nerve in the both directions so as to get a bisimplicial algebra,
then apply either the diagonal or the Artin-Mazur codiagonal functor (cf.
\cite{Artin}) to get to a simplicial algebra and take the Moore complex.
Consequently, he showed in \cite[Proposition 5.1]{Arvasi} that the Moore
complex of this simplicial algebra is isomorphic to the mapping cone
complex, in the sense of Loday, of the crossed square and this mapping cone
has a 2-crossed module structure of algebras.

Note that the construction given by Arvasi summarised above preserves the
homotopy modules. In fact, Ellis in \cite{ellis} proved that the homotopies
of the crossed square
\begin{equation*}
\xymatrix{ L \ar[d]_{\lambda'} \ar[r]^{\lambda} & M \ar[d]^{\mu} \\ N
\ar[r]_{\nu} & R }
\end{equation*}%
are the homologies of the complex
\begin{equation*}
\xymatrix{L\ar[rr]^-{(-\lambda,\lambda ^{\prime })}&&M\rtimes
N\ar[rr]^-{\mu+\nu}&&R\ar[r]&0.}
\end{equation*}

\section{\label{simp-square}Crossed Squares from Simplicial Algebras}

In 1991, Porter, \cite{porter2}, described a functor from the category of
simplicial groups to that of crossed $n$-cubes defined by Ellis and Steiner
in \cite{Sten}, based on ideas of Loday (cf. \cite{Loday}). Crossed $n$%
-cubes in algebraic settings such as commutative algebras, Jordan algebras,
Lie algebras have been defined by Ellis in \cite{elis}. The notion of
crossed $n$-cube of commutative algebras defined by Ellis is undoubtedly
important, however in this paper, we use only the case $n=2,$ that is for
crossed squares. Hence, we do not remind the general definition of a crossed
$n$-cube. From \cite{elis} we can say that, in low dimensions, a crossed
1-cube is the same as a crossed module and a crossed 2-cube is the same as a
crossed square.

In \cite{Arvasi}, Arvasi adapt that description to give an obvious analogue
of the functor given by Porter, \cite{porter2}, for the commutative algebra
case.

In fact, the following result is the 2-dimensional case of a general
construction of a crossed $n$-cube of algebras from a simplicial algebra
given by Arvasi, \cite{Arvasi}, analogue to that given by Porter, \cite%
{porter2}, in the group case.

Let $\mathbf{E}$ be a simplicial algebra. Then the following diagram
\begin{equation*}
\xymatrix{ NE_2/\partial_3 NE_3 \ar[d]_{\partial_2'} \ar[r]^-{\partial_2} &
NE_1 \ar[d]^{\mu} \\ \overline{NE}_1 \ar[r]_{\mu'} & E_1 }
\end{equation*}%
is the underlying square of a crossed square. The extra structure is given
as follows;\newline
$NE_{1}=\ker d_{0}^{1}$ and $\overline{NE}_{1}=\ker d_{1}^{1}$. Since $E_{1}$
acts on $NE_{2}/\partial _{3}NE_{3},\overline{NE}_{1}$ and $NE_{1}$, there
are actions of $\overline{NE}_{1}$ on $NE_{2}/\partial _{3}NE_{3}$ and $%
NE_{1}$ via $\mu ^{\prime }$, and $NE_{1}$ acts on $NE_{2}/\partial
_{3}NE_{3}$ and $\overline{NE}_{1}$ via $\mu $. Both $\mu $ and $\mu
^{\prime }$ are inclusions, and all actions are given by multiplication. The
$h$-map is
\begin{equation*}
\begin{tabular}{crcl}
$h$ : & $NE_{1}\times \overline{NE}_{1}$ & $\longrightarrow $ & $%
NE_{2}/\partial _{3}NE_{3}$ \\
& $(x,\overline{y})$ & $\longmapsto $ & $h(x,\overline{y}%
)=s_{1}x(s_{1}y-s_{0}y)+\partial _{3}NE_{3}$.%
\end{tabular}%
\end{equation*}%
Here $x$ and $y$ are in $NE_{1}$ as there is a natural bijection between $%
NE_{1}$ and $\overline{NE}_{1}.$ This is clearly functorial and we denote it
by;

\begin{center}
$\mathbf{M(}-\mathbf{,}2\mathbf{)}$\textbf{\ }$\mathbf{:SimpAlg}%
\longrightarrow \mathbf{Crs}$\textbf{$^{2}$}.
\end{center}

\section{\label{2crs-QM}Quadratic Modules from 2-Crossed Modules}

As we mentioned in introduction, Baues defined the quadratic module of
groups as an algebraic model of connected 3-types. In this section, we
define the commutative algebra version of this structure and we define a
functor from the category of 2-crossed modules to that of quadratic modules
of algebras. We should recall some basic information before giving the
definition of a quadratic module.

Recall that a \textit{pre-crossed module} is a homomorphism $\partial
:C\rightarrow R$ together with an action of $R$ on $C,$ written $c\cdot r$
for $r\in R$ and $c\in C,$ satisfying the condition $\partial (c\cdot
r)=\partial (c)r$ for all $r\in R$ and $c\in C.$

For an algebra $C,$ $C/C^{2}$ is the quotient of the algebra $C$ by its
ideal of squares. Then, there is a functor from the category of $\mathbf{k}$%
-algebras to the category of $\mathbf{k}$-modules. This functor goes from $C$
to $C/C^{2},$ plays the role of abelianization in the category of $\mathbf{k}
$-algebras. As modules are often called singular algebras (e.g., in the
theory of singular extensions) we shall call this functor \textquotedblleft
singularisation\textquotedblright .

Now, we generalise these notions to the pre-crossed modules. Let $\partial
:C\rightarrow R$ be a pre-crossed module and let $P_{1}(\partial )=C$ and
let $P_{2}(\partial )$ be the Peiffer ideal of $C$ generated by elements of
the form
\begin{equation*}
\left\langle x,y\right\rangle =xy-x\cdot \partial y
\end{equation*}%
which is called the\textit{\ Peiffer element} for $x,y\in C.$

A \textit{nil(2)-module} is a pre-crossed module $\partial :C\rightarrow R$
with an additional \textquotedblleft nilpotency\textquotedblright\
condition. This condition is $P_{3}(\partial )=0,$ where $P_{3}(\partial )$
is the ideal of $\ C$ generated by Peiffer elements $\left\langle
x_{1},x_{2},x_{3}\right\rangle $ of length $3.$

The homomorphism
\begin{equation*}
\partial ^{cr}:C^{cr}=C/P_{2}(\partial )\longrightarrow R
\end{equation*}%
is the crossed module associated to the pre-crossed module $\partial
:C\rightarrow R$, since
\begin{equation*}
\begin{array}{rll}
\left\langle [x],[y]\right\rangle & = (x+P_{2}(\partial ))(y+P_{2}(\partial
))-(x+P_{2}(\partial ))\cdot \partial ^{cr}(y+P_{2}(\partial )) &  \\
& = (xy+P_{2}(\partial ))-(x+P_{2}(\partial ))\cdot \partial
^{cr}(y+P_{2}(\partial )) &  \\
& = xy-x\cdot \partial (y)+P_{2}(\partial ) &  \\
& = P_{2}(\partial )\text{ \quad (}\because \text{ }\left\langle
x,y\right\rangle \in P_{2}(\partial )\text{ )} &  \\
& = [0] &
\end{array}%
\end{equation*}%
for $[x]=x+P_{2}(\partial ),$ $[y]=y+P_{2}(\partial )\in C^{cr}$.

Now, we can give the definition of a quadratic module of algebras.

\begin{defn}
A \textit{quadratic module} $(\omega ,\delta ,\partial )$ is a diagram
\begin{equation*}
\xymatrix{ & C\otimes C\ar[dl]_{\omega} \ar[d]^{w} \\L \ar[r]_{\delta} & M
\ar[r]_{\partial} &N }
\end{equation*}%
of homomorphisms of algebras such that the following axioms are satisfied.%
\newline
$\mathbf{QM1)}$- The homomorphism $\partial :M\rightarrow N$ is a
nil(2)-module and the quotient map $M\twoheadrightarrow C=$ $%
M_{{}}^{cr}/(M_{{}}^{cr})^{2}$ is given by $x\mapsto \lbrack x],$ where $%
[x]\in C$ denotes the class represented by $x\in M$. The map $w$ is defined
by Peiffer multiplication, i.e., $w([x]\otimes \lbrack y])=xy-x\cdot
\partial (y).$\newline
$\mathbf{QM2)}$- The homomorphisms $\delta $ and $\partial $ satisfy $\delta
\partial =0$ and the quadratic map $\omega $ is a lift of the map $w$, that
is $\delta \omega =w$ or equivalently
\begin{equation*}
\delta \omega ([x]\otimes \lbrack y])=w([x]\otimes \lbrack y])=xy-x\cdot
\partial (y)
\end{equation*}%
for $x,y\in M$.\newline
$\mathbf{QM3)}$- $L$ is a $N$-algebra and all homomorphisms of the diagram
are equivariant with respect to the action of $N$. Moreover, the action of $%
N $ on $L$ satisfies the following equality
\begin{equation*}
a\cdot \partial (x)=\omega ([\delta a]\otimes \lbrack x]+[x]\otimes \lbrack
\delta a])
\end{equation*}%
for $a\in L,x\in N$.\newline
$\mathbf{QM4)}$- For $a,b\in L$,
\begin{equation*}
\omega ([\delta a]\otimes \lbrack \delta b])=ab.
\end{equation*}
\end{defn}

A map $\varphi :(\omega ,\delta ,\partial )\rightarrow (\omega ^{\prime
},\delta ^{\prime },\partial ^{\prime })$ between quadratic modules is given
by a commutative diagram, $\varphi =(l,m,n)$
\begin{equation*}
\xymatrix{ C\otimes C\ar[d]_{\varphi_{\ast}\otimes\varphi_{\ast}}
\ar[r]^-{\omega} & L \ar[d]_{l} \ar[r]^{\delta} & M \ar[d]_{m}
\ar[r]^{\partial} & N\ar[d]_{n} \\ C^{\prime }\otimes C^{\prime }
\ar[r]_-{\omega^{\prime}} & L^{\prime} \ar[r]_{\delta^{\prime}} & M^{\prime}
\ar[r]_{\partial^{\prime}} & N^{\prime} }
\end{equation*}%
where $(m,n)$ is a map between pre-crossed modules which induces $\varphi
_{\ast }:C\rightarrow C^{\prime }$ and where $l$ is an $n$-equivariant
homomorphism. Let $\mathbf{QM}$ be the category of quadratic modules and of
maps as in above diagram.

Now, we construct a functor from the category of 2-crossed modules to the
category of quadratic modules.

Let
\begin{equation*}
\xymatrix{C_2 \ar[r]^-{\partial _{2}}&C_1 \ar[r]^-{\partial _{1}}&C_0}
\end{equation*}%
be a 2-crossed module of algebras. Let $P_{3}$ be the ideal of $C_{1}$
generated by elements of the form%
\begin{equation*}
\begin{array}{ccc}
\left\langle \left\langle x,y\right\rangle ,z\right\rangle & = &
\left\langle x,y\right\rangle z-\left\langle x,y\right\rangle \cdot \partial
_{1}z%
\end{array}%
\end{equation*}%
and%
\begin{equation*}
\begin{array}{ccc}
\left\langle x,\left\langle y,z\right\rangle \right\rangle & = &
x\left\langle y,z\right\rangle -x\cdot \partial _{1}(\left\langle
y,z\right\rangle )%
\end{array}%
\end{equation*}%
for $x,y,z\in C_{1}.$ Let
\begin{equation*}
q_{1}:C_{1}\rightarrow C_{1}/P_{3}
\end{equation*}%
be a quotient map and let $M=C_{1}/P_{3}$ be the quotient algebra. Since $%
\partial _{1}$ is a pre-crossed module, we obtain
\begin{equation*}
\begin{array}{lll}
\partial _{1}(\left\langle \left\langle x,y\right\rangle ,z\right\rangle ) &
= & \partial _{1}(\left\langle x,y\right\rangle z-\left\langle
x,y\right\rangle \cdot \partial _{1}z) \\
& = & \partial _{1}((xy-x\cdot \partial _{1}y)z-(xy-x\cdot \partial
_{1}y)\cdot \partial _{1}z) \\
& = & \partial _{1}x\partial _{1}y\partial _{1}z-\partial _{1}x\partial
_{1}y\partial _{1}z-\partial _{1}x\partial _{1}y\partial _{1}z+\partial
_{1}x\partial _{1}y\partial _{1}z \\
& = & 0%
\end{array}%
\end{equation*}%
and
\begin{equation*}
\begin{array}{lll}
\partial _{1}(\left\langle x,\left\langle y,z\right\rangle \right\rangle ) &
= & \partial _{1}(x\left\langle y,z\right\rangle -x\cdot \partial
_{1}\left\langle y,z\right\rangle ) \\
& = & \partial _{1}(x(yz-y\cdot \partial _{1}z)-x\cdot \partial
_{1}(yz-y\cdot \partial _{1}z)) \\
& = & \partial _{1}(x(yz-y\cdot \partial _{1}z)-x\cdot (\partial
_{1}y\partial _{1}z-\partial _{1}y\partial _{1}z)) \\
& = & \partial _{1}(x(yz-y\cdot \partial _{1}z)) \\
& = & \partial _{1}x\partial _{1}y\partial _{1}z-\partial _{1}x\partial
_{1}y\partial _{1}z \\
& = & 0%
\end{array}%
\end{equation*}%
for $x,y,z\in C_{1}$. That is, $\partial _{1}(P_{3})=0$. Thus, the map $%
\partial :M\rightarrow C_{0}$ given by $\partial (x+P_{3})=\partial _{1}(x),$
for all $x\in C_{1},$ is a well defined homomorphism since $\partial
_{1}(P_{3})=0$. Therefore, we can write a commutative diagram
\begin{equation*}
\xymatrix{ C_1 \ar[rr]^{\partial_1} \ar[dr]_{q_1} & & C_0 \\ & M
\ar[ur]_{\partial}}
\end{equation*}%
where $q_{1}:C_{1}\rightarrow M$ \ is the quotient map.

Let $P_{3}^{\prime }$ be the ideal of $C_{2}$ generated by elements of the
form
\begin{equation*}
\{x\otimes \left\langle y,z\right\rangle \}\text{ and }\{\left\langle
x,y\right\rangle \otimes z\},
\end{equation*}%
where $\{-\otimes -\}$ is the Peiffer lifting map. We have
\begin{equation*}
L=C_{2}/P_{3}^{\prime }.
\end{equation*}%
We can write from $\mathbf{2CM1)}$
\begin{equation*}
\partial _{2}\{x\otimes \left\langle y,z\right\rangle \}=\left\langle
x,\left\langle y,z\right\rangle \right\rangle
\end{equation*}%
and
\begin{equation*}
\partial _{2}\{\left\langle x,y\right\rangle \otimes z\}=\left\langle
\left\langle x,y\right\rangle ,z\right\rangle
\end{equation*}%
and we thus obtain $\partial _{2}(P_{3}^{\prime })=P_{3}$. Then, $\delta
:L\longrightarrow M$ given by $\delta (l+P_{3}^{\prime })=\partial
_{2}l+P_{3}$ is a well defined homomorphism. Indeed, if $l+P_{3}^{\prime
}=l^{\prime }+P_{3}^{\prime }$ then $l-l^{\prime }\in P_{3}^{\prime }$, and $%
\partial _{2}(l-l^{\prime })\in \partial _{2}(P_{3}^{\prime })$. Since $%
\partial _{2}(P_{3}^{\prime })=P_{3}$, we obtain $\partial _{2}(l)-\partial
_{2}(l^{\prime })\in P_{3}$, that is, $\partial _{2}l+P_{3}=\partial
_{2}l^{\prime }+P_{3}$.

Let
\begin{equation*}
C=\frac{M^{cr}}{(M^{cr})^{2}}.
\end{equation*}%
Thus we get the following commutative diagram;
\begin{equation*}
\xymatrix{ & C\otimes C\ar[dl]_{\omega} \ar[d]^{w} \\ L \ar[r]^{\delta } & M
\ar[r]^{\partial } & N \ar@{=}[d] \\ C_2 \ar[u]_{q_2} \ar[r]_{\partial_2} &
C_1 \ar[u]_{q_1} \ar[r]_{\partial_1} & C_0 }
\end{equation*}%
where $q_{1}$ and $q_{2}$ are the quotient maps. The quadratic map
\begin{equation*}
\omega :C\otimes C\longrightarrow L
\end{equation*}%
is given by the Peiffer lifting map, namely
\begin{equation*}
\omega \left( \lbrack q_{1}x]\otimes \lbrack q_{1}y]\right) =q_{2}\{x\otimes
y\}
\end{equation*}%
for all $x,y\in C_{1}$, $q_{1}x,q_{1}y\in M$ and $[q_{1}x]\otimes \lbrack
q_{1}y]\in C\otimes C.$

\begin{prop}
\label{quad} The diagram
\begin{equation*}
\xymatrix{ & C\otimes C\ar[dl]_{\omega} \ar[d]^{w} \\ L \ar[r]_{\delta} & M
\ar[r]_{\partial} & N }
\end{equation*}%
is a quadratic module of algebras.
\end{prop}

\begin{pf}
We show that all axioms of quadratic module are verified.

$\mathbf{QM1)}$\textbf{-} Since the triple Peiffer elements in $%
M=C_{1}/P_{3} $ are trivial, the $\partial :M\rightarrow N$ is a $nil(2)$%
-module. Indeed, for $x+P_{3},$ $y+P_{3},$ $z+P_{3}\in M,$
\begin{equation*}
\begin{array}{lll}
\left\langle x+P_{3},\left\langle y+P_{3},z+P_{3}\right\rangle \right\rangle
& = & \left\langle x,\left\langle y,z\right\rangle \right\rangle +P_{3} \\
& = & 0+P_{3}\quad \text{ (by }\left\langle x,\left\langle y,z\right\rangle
\right\rangle \in P_{3}\text{)} \\
& = & P_{3}%
\end{array}%
\end{equation*}%
and
\begin{equation*}
\begin{array}{lll}
\left\langle \left\langle x+P_{3},y+P_{3}\right\rangle ,z+P_{3}\right\rangle
& = & \left\langle \left\langle x,y\right\rangle ,z\right\rangle +P_{3} \\
& = & 0+P_{3}\quad \text{ (by }\left\langle \left\langle x,y\right\rangle
,z\right\rangle \in P_{3}\text{)} \\
& = & P_{3}.%
\end{array}%
\end{equation*}

$\mathbf{QM2)}$\textbf{-} For $q_{1}x,q_{1}y\in M$ and $[q_{1}x]\otimes
\lbrack q_{1}y]\in C\otimes C,$ we have%
\begin{equation*}
\begin{array}{ccl}
\delta \omega \left( \lbrack q_{1}x]\otimes \lbrack q_{1}y]\right)  & = &
\delta q_{2}\{x\otimes y\} \\
& = & q_{1}\partial _{2}\{x\otimes y\} \\
& = & q_{1}xq_{1}y-q_{1}x\cdot \partial (q_{1}y)\text{ \quad (by }\mathbf{%
2CM1)}\text{)} \\
& = & w([q_{1}x]\otimes \lbrack q_{1}y]).%
\end{array}%
\end{equation*}

$\mathbf{QM3)}$\textbf{-} For $q_{2}a\in L$ and $q_{1}x\in M$, and $[\delta
q_{2}a]\otimes \lbrack q_{1}x]\in C\otimes C,$ we have
\begin{equation*}
\begin{array}{ccl}
\omega \left( \lbrack \delta q_{2}a]\otimes \lbrack q_{1}x]\right)  & = &
\omega \left( \lbrack q_{1}\partial _{2}a]\otimes \lbrack q_{1}x]\right)  \\
& = & q_{2}\{\partial _{2}a\otimes x\} \\
& = & -q_{2}a\cdot x+q_{2}a\cdot \partial q_{1}(x)\text{\quad (by }\mathbf{%
2CM4)}\text{)}%
\end{array}%
\end{equation*}%
and%
\begin{equation*}
\begin{array}{ccl}
\omega \left( \lbrack q_{1}x]\otimes \lbrack \delta q_{2}a]\right)  & = &
\omega \left( \lbrack q_{1}x]\otimes \lbrack q_{1}\partial _{2}a]\right)  \\
& = & q_{2}\{x\otimes \partial _{2}a\} \\
& = & q_{2}a\cdot x\text{\quad (by }\mathbf{2CM4)}\text{).}%
\end{array}%
\end{equation*}%
Then, we have%
\begin{equation*}
\omega \left( \lbrack \delta q_{2}a]\otimes \lbrack q_{1}x]+[q_{1}x]\otimes
\lbrack \delta q_{2}a]\right) =q_{2}a\cdot \partial q_{1}(x).
\end{equation*}

$\mathbf{QM4)}$\textbf{-} For $q_{2}a,q_{2}b\in L$, and $[\delta
q_{2}a]\otimes \lbrack \delta q_{2}b]\in C\otimes C,$ we have%
\begin{equation*}
\begin{array}{ccl}
\omega \left( \lbrack \delta q_{2}a]\otimes \lbrack \delta q_{2}b]\right)  &
= & \omega \left( \lbrack q_{1}\partial _{2}a]\otimes \lbrack q_{1}\partial
_{2}b]\right)  \\
& = & q_{2}\{\partial _{2}a\otimes \partial _{2}b\} \\
& = & (q_{2}a)(q_{2}b)\text{\quad (by }\mathbf{2CM2)}\text{).}%
\end{array}%
\end{equation*}
\end{pf}

Therefore, we have defined a functor from the category of 2-crossed modules
to that of quadratic modules of algebras. We denote it by
\begin{equation*}
\Lambda :\mathbf{X}_{2}\mathbf{Mod\longrightarrow QM.}
\end{equation*}

Now, we show that the functor $\Lambda $ described above preserves the
homotopy modules.

\begin{prop}
\label{ho2} Let
\begin{equation*}
\xymatrix{C_2 \ar[r]^-{\partial _{2}}&C_1 \ar[r]^-{\partial _{1}}&C_0}
\end{equation*}%
be the 2-crossed module and $\pi _{i}$ be its homotopy modules for all $%
i\geqslant 0$. Let $\pi _{i}^{\prime }$ be the homotopy modules of its
associated quadratic module
\begin{equation*}
\xymatrix{ & C\otimes C\ar[dl]_{\omega} \ar[d]^{w} \\ L \ar[r]_{\delta } & M
\ar[r]_-{\partial} &N.}
\end{equation*}%
Then, $\pi _{i}\cong \pi _{i}^{\prime }$ for all $i\geqslant 0$.
\end{prop}

\begin{pf}
The homotopy modules of the 2-crossed module are
\begin{equation*}
\pi _{i}=\left\{
\begin{array}{lll}
C_{0}/\partial _{1}(C_{1}) &  & i=1, \\
\ker \partial _{1}/\mathrm{Im}\partial _{2} &  & i=2, \\
\ker \partial _{2} &  & i=3, \\
0 &  & i=0,\text{ }i>3%
\end{array}%
\right.
\end{equation*}%
and the homotopy modules of its associated quadratic module are
\begin{equation*}
\pi _{i}^{\prime }=\left\{
\begin{array}{lll}
N/\partial (M) &  & i=1, \\
\ker \partial /\mathrm{Im}\delta &  & i=2, \\
\ker \delta &  & i=3, \\
0 &  & i=0,\text{ }i>3.%
\end{array}%
\right.
\end{equation*}%
Now, we show that $\pi _{i}\cong \pi _{i}^{\prime }$ for all $i\geqslant 0$.
Since $\partial _{1}(P_{3})=0$, $M=C_{1}/P_{3}$, $N=C_{0}$ and $\partial
(M)\cong \partial _{1}(C_{1})$, clearly we have
\begin{equation*}
\pi _{1}=C_{0}/\partial _{1}(C_{1})\cong N/\partial (M)=\pi _{1}^{\prime }.
\end{equation*}%
Since $\ker \partial \cong \dfrac{\ker \partial _{1}}{P_{3}}$ and $\mathrm{Im%
}\delta \cong \dfrac{\mathrm{Im}\partial _{2}}{P_{3}}$ so that we have
\begin{equation*}
\pi _{2}^{\prime }=\dfrac{\ker \partial }{\mathrm{Im}\delta }\cong \dfrac{%
\ker \partial _{1}/P_{3}}{\mathrm{Im}\partial _{2}/P_{3}}\cong \dfrac{\ker
\partial _{1}}{\mathrm{Im}\partial _{2}}=\pi _{2}.
\end{equation*}%
Now, we show that $\pi _{3}\cong \pi _{3}^{\prime }$. Consider that
\begin{equation*}
\pi _{3}^{\prime }=\{x+P_{3}^{\prime }:\partial _{2}(x)\in P_{3}\}.
\end{equation*}%
For an element $x+P_{3}^{\prime }$ of $\pi _{3}^{\prime }$, we show that
there is an element $x^{\prime }+P_{3}^{\prime }$ of $\pi _{3}^{\prime }$
such that $x+P_{3}^{\prime }=x^{\prime }+P_{3}^{\prime }$ and $x^{\prime
}\in \ker \partial _{2}$. In fact, we observe from $\mathbf{2CM1)}$ that $%
\partial _{2}\{\left\langle x,y\right\rangle \otimes z\}=\left\langle
\left\langle x,y\right\rangle ,z\right\rangle $ and $\partial _{2}\{x\otimes
\left\langle y,z\right\rangle \}=\left\langle x,\left\langle
y,z\right\rangle \right\rangle $, and we have $\partial _{2}(P_{3}^{\prime
})=P_{3}$. Hence $\partial _{2}(x)\in P_{3}$ implies $\partial
_{2}(x)=\partial _{2}(w),$ $w\in P_{3}^{\prime };$ thus $\partial
_{2}(x-w)=0;$ then take $x^{\prime }=x-w,$ so that $x+P_{3}^{\prime
}=x^{\prime }+P_{3}^{\prime }$ and $\partial _{2}(x^{\prime })=0.$ Define $%
\mu :\pi _{3}^{\prime }\rightarrow \pi _{3}$ by $(x^{\prime }+P_{3}^{\prime
})\mapsto x^{\prime }$ and $\nu :\pi _{3}\rightarrow \pi _{3}^{\prime }$ by $%
x\mapsto x+P_{3}$. Then $\mu $ and $\nu $ are inverse bijections, that is,
we have $\pi _{3}\cong \pi _{3}^{\prime }$. Thus, the homotopy modules of
the 2-crossed module are isomorphic to that of its associated quadratic
module.
\end{pf}

\section{\label{simp-QM}Quadratic Modules from Simplicial Algebras}

Baues, \cite{baus2}, defined a functor from simplicial groups to quadratic
modules. In the construction of this functor, to define the quadratic map $%
\omega $, Baues in \cite{baus2} used the Conduch\'{e}'s Peiffer lifting map $%
\{-,-\}$ given for the construction of a 2-crossed module from a simplicial
group (cf.  \cite{con}).

In this section, we construct a functor from the category of simplicial
algebras to the category of quadratic modules in terms of hypercrossed
complex pairings analogously to that given by Baues in \cite{baus2}. In the
construction, to define the quadratic map $\omega $, we use the Peiffer
lifting map given by Arvasi and Porter in \cite[Proposition 5.2]{Arvasi1}.

Now, we will construct a functor from simplicial algebras to quadratic
modules by using the $C_{\alpha, \beta}$ functions. We will use the $%
C_{\alpha, \beta}$ functions in verifying the axioms of quadratic module.

Let $\mathbf{E}$ be a simplicial algebra with Moore complex $\mathbf{NE}$.
We will obtain a quadratic module by using the following commutative diagram
\begin{equation*}
\xymatrix{ & NE_1\times NE_1\ar[dl]_{\omega'} \ar[d]^{w} \\ NE_2/\partial_3
(NE_3 )\ar[r]_-{\overline{\partial_2}} & NE_1 \ar[r]_-{\partial_1} & NE_0}
\end{equation*}%
where the map $w$ is given by
\begin{equation*}
w(x,y)=xy-x\cdot \partial _{1}y
\end{equation*}%
for $x,y\in NE_{1}$ and the map $\omega ^{\prime }$ is given by
\begin{equation*}
\omega ^{\prime }(x,y)=\overline{s_{1}x(s_{1}y-s_{0}y)}%
=s_{1}x(s_{1}y-s_{0}y)+\partial _{3}(NE_{3})
\end{equation*}%
for $x,y\in NE_{1}.$

Let $P_{3}(\partial _{1})$ be the ideal of $NE_{1}$ generated by elements of
the form
\begin{equation*}
\left\langle x,\left\langle y,z\right\rangle \right\rangle \text{ and }%
\left\langle \left\langle x,y\right\rangle ,z\right\rangle
\end{equation*}%
for $x,y,z\in NE_{1}$. Since $\partial _{1}$ is a pre-crossed module, we can
write $\partial _{1}(P_{3}(\partial _{1}))=0.$ Then, $\partial
:NE_{1}/P_{3}(\partial _{1})\longrightarrow NE_{0}$ given by $\partial
(x+P_{3}(\partial _{1}))=\partial _{1}x$ for $x\in NE_{1}$ is a well defined
homomorphism. Thus, we obtain the following commutative diagram
\begin{equation*}
\xymatrix{NE_1\ar[rr]^{\partial_1}\ar[dr]_{q_1}&&NE_0\\
&NE_1/P_3(\partial_1)\ar[ur]_{\partial}&}
\end{equation*}%
where $q_{1}$ is the quotient map.

Let $P_{3}^{\prime }(\partial _{1})$ be an ideal of $NE_{2}/\partial
_{3}(NE_{3})$ generated by the formal Peiffer elements of $x,y,z$, i.e.,
generated by elements of the form
\begin{equation*}
\omega ^{\prime }(\left\langle x,y\right\rangle ,z)=s_{1}(\left\langle
x,y\right\rangle )(s_{1}z-s_{0}z)
\end{equation*}%
and
\begin{equation*}
\omega ^{\prime }(x,\left\langle y,z\right\rangle
)=s_{1}(x)(s_{1}(\left\langle y,z\right\rangle )-s_{0}(\left\langle
y,z\right\rangle ))
\end{equation*}%
for $x,y,z\in NE_{1}.$ Notice that
\begin{equation*}
\begin{array}{lll}
\overline{\partial _{2}}\omega ^{\prime }(x,\left\langle y,z\right\rangle )
& = & d_{2}(s_{1}(x)(s_{1}(\left\langle y,z\right\rangle
)-s_{0}(\left\langle y,z\right\rangle ))) \\
& = & x(\left\langle y,z\right\rangle )-xs_{0}d_{1}(\left\langle
y,z\right\rangle ) \\
& = & \left\langle x,\left\langle y,z\right\rangle \right\rangle%
\end{array}%
\end{equation*}%
and
\begin{equation*}
\begin{array}{lll}
\overline{\partial _{2}}\omega ^{\prime }(\left\langle x,y\right\rangle ,z)
& = & d_{2}(s_{1}(\left\langle x,y\right\rangle )(s_{1}(z)-s_{0}(z))) \\
& = & \left\langle x,y\right\rangle (z-s_{0}d_{1}z) \\
& = & \left\langle \left\langle x,y\right\rangle ,z\right\rangle .%
\end{array}%
\end{equation*}%
We obtain $\overline{\partial _{2}}(P_{3}^{\prime }(\partial
_{1}))=P_{3}(\partial _{1}).$ Let
\begin{equation*}
M=NE_{1}/P_{3}(\partial _{1})
\end{equation*}%
and
\begin{equation*}
L=(NE_{2}/\partial _{3}(NE_{3}))/P_{3}^{\prime }(\partial _{1})
\end{equation*}%
be the quotient algebras. We thus see that the map
\begin{equation*}
\delta :(NE_{2}/\partial _{3}NE_{3}))/P_{3}^{\prime }(\partial
_{1})\longrightarrow NE_{1}/P_{3}(\partial _{1})
\end{equation*}%
given by $\delta (a+P_{3}^{\prime }(\partial _{1}))=\overline{\partial _{2}}%
(a)+P_{3}(\partial _{1})$ is a well defined homomorphism since $\overline{%
\partial _{2}}(P_{3}^{\prime }(\partial _{1}))=P_{3}(\partial _{1})$.

We obtain the following commutative diagram,
\begin{equation*}
\xymatrix{ & C\otimes C\ar[dl]_{\omega} \ar[d]^{w} \\ L \ar[r]_-{\delta } &
M \ar[r]_-{\partial} & N \ar@{=}[d] \\ NE_2/\partial_3( NE_3 ) \ar[u]^{q_2}
\ar[r]_-{\overline{\partial_2}} & NE_1 \ar[u]_{q_1} \ar[r]_-{\partial_1} &
NE_0 }
\end{equation*}%
where $C=M^{cr}/(M^{cr})^{2}$ and $q_{1},q_{2}$ are the quotient maps and
the quadratic map $\omega $ can be given by
\begin{equation*}
\omega ([q_{1}x]\otimes \lbrack q_{1}y])=q_{2}\omega ^{\prime
}(x,y)=q_{2}\left( s_{1}x(s_{1}y-s_{0}y)+\partial _{3}(NE_{3})\right)
\end{equation*}%
for $q_{1}x,q_{1}y\in M$ and $[q_{1}x]\otimes \lbrack q_{1}y]\in C\otimes C$.

\begin{prop}
\label{simp}The diagram
\begin{equation*}
\xymatrix{ & C\otimes C\ar[dl]_{\omega} \ar[d]^{w} \\ L \ar[r]_-{\delta } &
M \ar[r]_-{\partial} & N }
\end{equation*}%
is a quadratic module of algebras.
\end{prop}

\begin{pf}
We show that all axioms of quadratic module are verified by using the $%
C_{\alpha ,\beta }$ functions. We display the elements omitting the
overlines in our calculation to save complication.

$\mathbf{QM1)}$- Clearly, $\partial :M\rightarrow N$ is a nil(2)-module.
Because, for $x+P_{3}(\partial _{1}),$ $y+P_{3}(\partial _{1}),$ $%
z+P_{3}(\partial _{1})\in M=NE_{1}/P_{3}(\partial _{1}),$
\begin{equation*}
\begin{array}{lll}
\left\langle x+P_{3}(\partial _{1}),\left\langle y+P_{3}(\partial
_{1}),z+P_{3}(\partial _{1})\right\rangle \right\rangle & = & \left\langle
x,\left\langle y,z\right\rangle \right\rangle +P_{3}(\partial _{1}) \\
& = & 0+P_{3}(\partial _{1})\text{ (by }\left\langle x,\left\langle
y,z\right\rangle \right\rangle \in P_{3}(\partial _{1})\text{)} \\
& = & P_{3}(\partial _{1})%
\end{array}%
\end{equation*}%
and similarly
\begin{equation*}
\begin{array}{lll}
\left\langle \left\langle x+P_{3}(\partial _{1}),y+P_{3}(\partial
_{1})\right\rangle ,z+P_{3}(\partial _{1})\right\rangle & = & \left\langle
\left\langle x,y\right\rangle ,z\right\rangle +P_{3}(\partial _{1}) \\
& = & 0+P_{3}(\partial _{1})\text{ (}\because \text{ }\left\langle
x,\left\langle y,z\right\rangle \right\rangle \in P_{3}(\partial _{1})\text{)%
} \\
& = & P_{3}(\partial _{1}).%
\end{array}%
\end{equation*}

$\mathbf{QM2)}$- For $x,y\in NE_{1},$ $q_{1}x,q_{1}y\in M$ and $%
[q_{1}x]\otimes \lbrack q_{1}y]\in C\otimes C,$ we have
\begin{equation*}
\begin{array}{lll}
\delta \omega ([q_{1}x]\otimes \lbrack q_{1}y]) & = & \delta q_{2}\left(
s_{1}x(s_{1}y-s_{0}y)\right)  \\
& = & q_{1}(d_{2}(s_{1}x(s_{1}y-s_{0}y)))\quad \text{(}\because \text{ }%
\delta q_{2}=q_{1}\overline{\partial _{2}}) \\
& = & q_{1}(x)q_{1}(y)-q_{1}(x)\cdot \overline{\partial _{1}}(q_{1}y)\quad
\text{(}\because \text{ }\partial q_{1}=d_{1}) \\
& = & w([q_{1}x]\otimes \lbrack q_{1}y]).%
\end{array}%
\end{equation*}

$\mathbf{QM3)}$- For $q_{2}a\in L$ and $q_{1}x\in M,$ we have
\begin{equation*}
\omega \left( \lbrack \delta q_{2}a]\otimes \lbrack q_{1}x]\right) =\omega
\left( \lbrack q_{1}\overline{\partial _{2}}a]\otimes \lbrack q_{1}x]\right)
\end{equation*}%
and
\begin{equation*}
\omega \left( \lbrack q_{1}\overline{\partial _{2}}a]\otimes \lbrack
q_{1}x]\right) =q_{2}(s_{1}d_{2}a(s_{1}x-s_{0}x)).
\end{equation*}%
On the other hand, from
\begin{equation*}
\partial _{3}(C_{(2,0)(1)}(x\otimes
a))=(s_{0}x-s_{1}x)s_{1}d_{2}a-(s_{0}x-s_{1}x)a\in \partial _{3}(NE_{3})
\end{equation*}%
we have
\begin{equation*}
\begin{array}{ccl}
\omega \left( \lbrack q_{1}\overline{\partial _{2}}a]\otimes \lbrack
q_{1}x]\right)  & \equiv  & s_{0}d_{1}x(q_{2}a)-s_{1}x(q_{2}a)\quad \func{mod%
}(\partial _{3}(NE_{3})) \\
& = & (q_{2}a)\cdot \overline{\partial _{1}}q_{1}x-x\cdot (q_{2}a)\quad
\text{(}\because \text{ }\overline{\partial _{1}}q_{1}=d_{1}).%
\end{array}%
\end{equation*}%
Similarly, we have
\begin{equation*}
\begin{array}{lll}
\omega \left( \lbrack q_{1}x]\otimes \lbrack \delta q_{2}a]\right)  & = &
\omega \left( \lbrack q_{1}x]\otimes \lbrack q_{1}\overline{\partial _{2}}%
a]\right)  \\
& = & q_{2}(s_{1}x(s_{1}d_{2}a-s_{0}d_{2}a)).%
\end{array}%
\end{equation*}%
From
\begin{equation*}
\partial _{3}(C_{(2,1)(0)}(x\otimes
a))=s_{1}x(s_{0}d_{2}a-s_{1}d_{2}a)+s_{1}xa\in \partial _{3}(NE_{3})
\end{equation*}%
we can write
\begin{equation*}
\begin{array}{ccl}
\omega \left( \lbrack q_{1}x]\otimes \lbrack q_{1}\overline{\partial _{2}}%
a]\right)  & \equiv  & x\cdot (q_{2}a)\quad \func{mod}(\partial
_{3}(NE_{3})).%
\end{array}%
\end{equation*}%
Therefore,
\begin{equation*}
\omega ([\delta q_{2}a]\otimes \lbrack q_{1}x]+[q_{1}x]\otimes \lbrack
\delta q_{2}a])=(q_{2}a)\cdot \partial (q_{1}x).
\end{equation*}

$\mathbf{QM4)}$- For $q_{2}a,q_{2}b\in L$ , we can write
\begin{equation*}
\omega \left( \lbrack \delta q_{2}a]\otimes \lbrack \delta q_{2}b]\right)
=q_{2}(s_{1}d_{2}a(s_{1}d_{2}b-s_{0}d_{2}b)).
\end{equation*}%
From
\begin{equation*}
\partial _{3}(C_{(1)(0)}(a\otimes
b))=s_{1}d_{2}a(s_{0}d_{2}b-s_{1}d_{2}b)+ab\in \partial _{3}(NE_{3})
\end{equation*}%
we obtain
\begin{equation*}
\begin{array}{ccl}
\omega \left( \lbrack q_{1}\overline{\partial _{2}}a]\otimes \lbrack q_{1}%
\overline{\partial _{2}}b]\right) & \equiv & q_{2}(a)q_{2}(b)\quad \func{mod}%
(\partial _{3}(NE_{3})).%
\end{array}%
\end{equation*}
\end{pf}

Therefore, we defined a functor from the category of simplicial algebras to
that of quadratic modules of algebras. We denote it by
\begin{equation*}
\Delta :\mathbf{SimpAlg}\longrightarrow \mathbf{QM}.
\end{equation*}

Alternatively, this proposition can be reproved differently, by making use
of the 2-crossed module constructed from a simplicial algebra by Arvasi and
Porter (cf. \cite{Arvasi1}). We now give a sketch of the argument. In \cite[%
Proposition 5.2]{Arvasi1}, it is shown that given a simplicial algebra $%
\mathbf{E}$, one can construct a 2-crossed module

\begin{equation}
\xymatrix{ NE_2/\partial_3(NE_3 \cap D_3) \ar[r]^-{\overline{\partial_2}} &
NE_1 \ar[r]^{\partial_1 } & NE_0 }  \tag{1}  \label{1}
\end{equation}%
where $\{x\otimes y\}=s_{1}x(s_{1}y-s_{0}y)+\partial _{3}(NE_{3}\cap D_{3})$
for $x,y\in NE_{1}.$ \newline
Clearly we have a commutative diagram
\begin{equation*}
\xymatrix{ NE_2/\partial_3(NE_3 \cap D_3)\ar[d]^{j}
\ar[r]^-{\overline{\partial_2}} & NE_1 \ar[r]^{\partial_1 }\ar@{=}[d] & NE_0
\ar@{=}[d]\\ NE_2/\partial_3(NE_3)\ar[r]&NE_1 \ar[r]^{\partial_1}&NE_0.}
\end{equation*}%
Consider now the quadratic module associated to the 2-crossed module (\ref{1}%
), as in Section \ref{2crs-QM} of this paper.
\begin{equation*}
\xymatrix{ & C\otimes C\ar[dl]_{\omega'} \ar[d] \\ L '\ar[r]^{\delta' } & M
\ar[r] & N \ar@{=}[d] \\ \dfrac{{NE_2}}{{\partial_3 (NE_3\cap D_3)}}
\ar[u]_{q_2} \ar[r]_-{\overline{\partial_2}} & NE_1 \ar[u]
\ar[r]_{\partial_1} & NE_0 }
\end{equation*}%
Then one can see that $L^{\prime }=A /\partial _{3}(NE_{3}\cap
D_{3}),$ where $A $ is the ideal of $NE_{2}$ generated by elements
of the form
\begin{equation*}
s_{1}(\left\langle x,y\right\rangle )(s_{1}z-s_{0}z)\text{ and }%
s_{1}(x)(s_{1}(\left\langle y,z\right\rangle )-s_{0}(\left\langle
y,z\right\rangle )).
\end{equation*}%
On the other hand we have, from Section \ref{2crs-QM}, $L=A
/\partial _{3}(NE_{3})$. Hence there is a map $i:L^{\prime
}\rightarrow L$ with
\begin{equation}
\omega =i\omega ^{\prime },\quad \delta ^{\prime }=\delta i.  \tag{2}
\label{2}
\end{equation}%
Since
\begin{equation*}
\xymatrix{ C\otimes C \ar[r]^-{\omega'}&L'\ar[r]^{\delta'} &M
\ar[r]^{\partial } & N }
\end{equation*}%
is, by construction a quadratic module, it is straightforward to check,
using (\ref{2}), that
\begin{equation*}
\xymatrix{ C\otimes C \ar[r]^-{\omega}&L\ar[r]^{\delta} &M \ar[r]^{\partial
} & N }
\end{equation*}%
is also a quadratic module.

Now we show that the functor $\Delta $ described above preserves the
homotopy modules.

\begin{prop}
Let $\mathbf{E}$ be a simplicial  algebra, let $\pi _{i}$ be the
homotopy modules of the classifying space of $\mathbf{E}$ and let
$\pi _{i}^{\prime }$ be the homotopy modules of its associated
quadratic module; then $\pi _{i}\cong \pi _{i}^{\prime }$ for
$i=0,1,2,3$.
\end{prop}

\begin{pf}
Let $\mathbf{E}$ be a simplicial algebra. The $n$th homotopy module of $%
\mathbf{E}$ is isomorphic to the $n$th homology of the Moore complex of $%
\mathbf{E}$, i.e., $\pi _{n}(\mathbf{E})\cong H_{n}(\mathbf{NE}).$Thus the
homotopy modules $\pi _{n}(\mathbf{E})=\pi _{n}$ of $\mathbf{E}$ are%
\begin{equation*}
\pi _{n}=\left\{
\begin{array}{lll}
NE_{0}/d_{1}(NE_{1}) &  & n=1, \\
&  &  \\
\dfrac{{\ker d_{1}\cap NE_{1}}}{{d_{2}(NE_{2})}} &  & n=2, \\
&  &  \\
\dfrac{\ker d_{2}\cap NE_{2}}{d_{3}({NE_{3}})} &  & n=3, \\
&  &  \\
0 &  & n=0, \text{or}\ n>3%
\end{array}%
\right.
\end{equation*}%
and the homotopy modules $\pi _{n}^{\prime }$ of its associated quadratic
module are
\begin{equation*}
\pi _{n}^{\prime }=\left\{
\begin{array}{lll}
NE_{0}/\partial (M) &  & n=1, \\
\ker \partial /\mathrm{\func{Im}}\delta  &  & n=2, \\
\ker \delta  &  & n=3, \\
0 &  & n=0, \text{or}\  n>3.%
\end{array}%
\right.
\end{equation*}%
We claim that $\pi _{n}^{\prime }\cong \pi _{n}$ for $n=1,2,3$. Since $%
M=NE_{1}/P_{3}(\partial _{1})$ and $\partial _{1}(P_{3}(\partial _{1}))=0$,
we have
\begin{equation*}
\partial (M)=\partial (NE_{1}/P_{3}(\partial _{1}))=d_{1}(NE_{1})
\end{equation*}%
and then
\begin{equation*}
\pi _{1}^{\prime }=NE_{0}/\partial (M)\cong NE_{0}/d_{1}(NE_{1})=\pi _{1}.
\end{equation*}%
Also $\ker \partial =\dfrac{\ker d_{1}\cap NE_{1}}{P_{3}(\partial _{1})}$
and $\mathrm{\func{Im}}\delta =d_{2}(NE_{2})/P_{3}(\partial _{1})$ so that
we have
\begin{equation*}
\pi _{2}^{\prime }=\frac{\ker \partial }{\func{Im}\delta }=\frac{(\ker
d_{1}\cap NE_{1})/P_{3}(\partial _{1})}{d_{2}(NE_{2})/P_{3}(\partial _{1})}%
\cong \frac{\ker d_{1}\cap NE_{1}}{d_{2}(NE_{2})}=\pi _{2}.
\end{equation*}%
The isomorphism between $\pi _{3}^{\prime }$ and $\pi _{3}$ can be proved
similarly to the proof of Proposition \ref{ho2}.
\end{pf}

\textbf{Remark: }Note that in the previous section, we have defined the
functor $\Lambda $ from the category of 2-crossed modules to that of
quadratic modules. As mentioned in Section \ref{simp-2crs} and above, Arvasi
and Porter, \cite{Arvasi1}, constructed an equivalence between simplicial
algebras with Moore complex of length 2 and 2-crossed modules of algebras.
We can summarise these statements in the following diagram
\begin{equation*}
\xymatrix{\textbf{X$_2$Mod}\ar@<-0.9ex>[dr]_{\Lambda}\ar[rr]^{%
\cite{Arvasi1}}&&\textbf{SimpAlg$_{\leq
2}$}\ar[ll]\ar@<0.9ex>[dl]^{\Delta}\\ &\textbf{QM}&}
\end{equation*}

\section{\label{square-QM}Quadratic Modules from Crossed Squares}

In this section we will define a functor from crossed squares to quadratic
modules of algebras. Our construction can be briefly explained as:

Given a crossed square of algebras, we consider the associated 2-crossed
module from Section \ref{square-2crs} (cf. \cite{Arvasi}), and then we build
the quadratic module corresponding to this 2-crossed module as in Section %
\ref{2crs-QM}. In other words, we are just composing two functors. In
particular, the homotopy type is clearly preserved, as it is preserved at
each step.

Now, recall that the first author in \cite{Arvasi} constructed a 2-crossed
module from a crossed square of commutative algebras
\begin{equation*}
\xymatrix{ L \ar[d]_{\lambda'} \ar[r]^{\lambda} & M \ar[d]^{\mu} \\ N
\ar[r]_{\nu} & R }
\end{equation*}%
as
\begin{equation}
\xymatrix{L\ar[rr]^-{(-\lambda,\lambda ^{\prime })}&&M\rtimes
N\ar[rr]^-{\mu+\nu}&&R }  \tag{3}  \label{(3)}
\end{equation}%
analogue to that given by Conduch\'{e} in the group case (cf. \cite{Con1}),
as we mention in Section \ref{square-2crs}.

Now let
\begin{equation*}
\begin{array}{c}
\xymatrix{ L\ar[r]^{\lambda}\ar[d]_{\lambda'}& M\ar[d]^{\mu}\\ N
\ar[r]_{\nu}&R}%
\end{array}%
\end{equation*}%
be a crossed square of algebras. Consider its associated 2-crossed module (%
\ref{(3)}). From this 2-crossed module, we can construct a quadratic module
as in Section \ref{2crs-QM},
\begin{equation*}
\xymatrix{ & C\otimes C\ar[dl]_{\omega} \ar[d]^{w} \\ L' \ar[r]_-{\delta} &
M' \ar[r]_-{\partial} & N}
\end{equation*}%
where $N=R,$ $M^{\prime }=(M\rtimes N)/P_{3},$ $L^{\prime }=L/P_{3}^{\prime
} $ , $C=(M^{\prime })^{cr}/((M^{\prime })^{cr})^{2}$.

The Peiffer elements in $M\rtimes N$ are given by
\begin{equation*}
\begin{array}{lll}
\left\langle (m,n),(c,a)\right\rangle & = & (m,n)(c,a)-(m,n)\cdot (\mu
(c)+\nu (a)) \\
& = & (mc+m\cdot \nu (a)+c\cdot \nu (n),na) \\
&  & -(m\cdot \mu (c)+m\cdot \nu (a),n\cdot \mu (c)+n\cdot \nu (a)) \\
& = & (c\cdot \nu (n),-n\cdot \mu (c)).%
\end{array}%
\end{equation*}%
We know from Section \ref{2crs-QM} that the ideal $P_{3}$ of $M\rtimes N$ is
generated by elements of the form%
\begin{equation*}
\left\langle \left\langle (m,n),(c,a)\right\rangle ,(m^{\prime },n^{\prime
})\right\rangle =(-m^{\prime }\cdot \nu (n\cdot \mu (c)),\mu (m^{\prime
})\cdot (n\cdot \mu (c)))
\end{equation*}%
and
\begin{equation*}
\left\langle (m,n),\left\langle (c,a),(m^{\prime },n^{\prime })\right\rangle
\right\rangle =((m^{\prime }\cdot \nu (a))\cdot \nu (n),-n\cdot \mu
(m^{\prime }\cdot \nu (a)))
\end{equation*}%
for $(m,n),(c,a),(m^{\prime },n^{\prime })\in M\rtimes N$. We thus have $%
(\mu +\nu )(P_{3})=0$ and the map $\partial :(M\rtimes N)/P_{3}\rightarrow R
$ is given by $\partial ((m,n)+P_{3})=\mu (m)+\nu (n).$

We also know that the ideal $P_{3}^{\prime }$ of $L$ is generated by
elements of the form
\begin{equation*}
h(m^{\prime },-n\cdot \mu (c))\text{ and }h(m^{\prime }\cdot \nu (a),n)
\end{equation*}%
for all $m^{\prime },c\in M$ and $n,a\in N.$ We have
\begin{equation*}
\begin{array}{lll}
(-\lambda ,\lambda ^{\prime })(h(m^{\prime },-n\cdot \mu (c))) & = &
(-m^{\prime }\cdot \nu (n\cdot \mu (c)),\mu (m^{\prime })\cdot (n\cdot \mu
(c))) \\
& = & \left\langle \left\langle (m,n),(c,a)\right\rangle ,(m^{\prime
},n^{\prime })\right\rangle%
\end{array}%
\end{equation*}%
and
\begin{equation*}
\begin{array}{lll}
(-\lambda ,\lambda ^{\prime })(h(m^{\prime }\cdot \nu (a),n)) & = &
((m^{\prime }\cdot \nu (a))\cdot \nu (n),-n\cdot \mu (m^{\prime }\cdot \nu
(a))) \\
& = & \left\langle (m,n),\left\langle (c,a),(m^{\prime },n^{\prime
})\right\rangle \right\rangle%
\end{array}%
\end{equation*}%
by the crossed square axioms. Thus the map $\delta :L^{\prime }\rightarrow
M^{\prime }$ is given by $\delta \left( l+P_{3}^{\prime }\right) =(-\lambda
l,\lambda ^{\prime }l)+P_{3}$ for all $l\in L$. The quadratic map
\begin{equation*}
\omega :C\otimes C\longrightarrow L^{\prime }
\end{equation*}%
can be given by
\begin{equation*}
\omega \left( \lbrack q_{1}(m,n)]\otimes \lbrack q_{1}(c,a)]\right)
=q_{2}(h(c,na))
\end{equation*}%
for all $(m,n),(c,a)$ $\in M\rtimes N$ , $q_{1}(m,n),$ $q_{1}(c,a)\in
M^{\prime }$ and $[q_{1}(m,n)]\otimes \lbrack q_{1}(c,a)]\in C\otimes C,$
and where $h$ is the $h$-map of the crossed square.

\begin{prop}
\label{crs} The diagram
\begin{equation*}
\xymatrix{ & C\otimes C\ar[dl]_{\omega} \ar[d]^{w} \\ L' \ar[r]_-{\delta} &
M' \ar[r]_-{\partial} & N }
\end{equation*}%
is a quadratic module of algebras.
\end{prop}

\begin{pf}
We show that all axioms of quadratic module are verified.

$\mathbf{QM1)}$\textbf{-} Obviously $\partial :M^{\prime }\rightarrow N$ is
a $nil(2)$-module, since the triple Peiffer elements in $M^{\prime }$ are
trivial. It can be proved similarly to the proof of Proposition \ref{quad}.

$\mathbf{QM2)}$\textbf{-} For $q_{1}(m,n)$ and $q_{1}(c,a)\in M^{\prime }$
and $[q_{1}(m,n)],[q_{1}(c,a)]\in C$
\begin{align*}
\delta \omega \left( \lbrack q_{1}(m,n)]\otimes \lbrack q_{1}(c,a)]\right) &
=(-\lambda ,\lambda ^{\prime })q_{2}(h(c,na)) \\
& =q_{1}((-\lambda h(c,na),\lambda ^{\prime }h(c,na)) \\
& =w([q_{1}(m,n)]\otimes \lbrack q_{1}(c,a)]).
\end{align*}

$\mathbf{QM3)}$\textbf{-} For $q_{1}(c,a),\delta (q_{2}l)\in M^{\prime }$
and $[\delta q_{2}l],[q_{1}(c,a)]\in C,$%
\begin{align*}
\omega ([\delta q_{2}l]\otimes \lbrack q_{1}(c,a)])+\omega
([q_{1}(c,a)]\otimes \lbrack \delta q_{2}l])& =\omega \lbrack q_{1}(-\lambda
l,\lambda ^{\prime }l)]\otimes \lbrack q_{1}(c,a)]+\omega \lbrack
q_{1}(c,a)]\otimes \lbrack q_{1}(-\lambda l,\lambda ^{\prime }l)] \\
& =q_{2}(h(c,(\lambda ^{\prime }l)a))+q_{2}(-\lambda l,a\lambda ^{\prime }l)
\\
& =q_{2}(h(c,\lambda ^{\prime }(l\cdot a))+q_{2}(h(-\lambda l,\lambda
^{\prime }(l\cdot a)))\text{ } \\
& =(\mu (c)+\nu (a))\cdot q_{2}(l).
\end{align*}

$\mathbf{QM4)}$\textbf{-} For $\delta (q_{2}l)=q_{1}(-\lambda l,\lambda
^{\prime }l)$ and $\delta (q_{2}l_{0})=q_{1}(-\lambda l_{0},\lambda ^{\prime
}l_{0})\in M^{\prime }$ and $[\delta q_{2}l],[\delta q_{2}l_{0}]\in C,$
\begin{align*}
\omega \left( \lbrack \delta q_{2}l]\otimes \lbrack \delta
q_{2}l_{0}]\right) & =\omega \lbrack q_{1}(-\lambda l,\lambda ^{\prime
}l)]\otimes \lbrack q_{1}(-\lambda l_{0},\lambda ^{\prime }l_{0})] \\
& =q_{2}(h(\lambda l,\lambda ^{\prime }(ll_{0}))) \\
& =(q_{2}l)(q_{2}l_{0}).
\end{align*}
\end{pf}

Thus, we would have defined a functor from the category of crossed squares
to that of quadratic modules of commutative algebras and we denote it by
\begin{equation*}
\Psi :\mathbf{Crs}^{2}\longrightarrow \mathbf{QM.}
\end{equation*}%
Thus, the diagram given in introduction was completed.

\begin{equation*}
\begin{array}{llllll}
\text{Z. Arvasi} &  &  &  &  & \text{E. Ulualan} \\
\text{Eski\c{s}ehir Osmangazi University} &  &  &  &  & \text{Dumlup\i nar
University} \\
\text{Science Faculty} &  &  &  &  & \text{Science Faculty} \\
\text{Department of Mathematics} &  &  &  &  & \text{Department of
Mathematics} \\
\text{26480, Eski\c{s}ehir-TURKEY} &  &  &  &  & \text{K\"{u}tahya-TURKEY}
\\
\text{e-mail: zarvasi@ogu.edu.tr} &  &  &  &  & \text{%
eulualan@dumlupinar.edu.tr}%
\end{array}%
\end{equation*}

\bigskip

\end{document}